\documentclass[12pt]{article}
\usepackage{amsmath,amsfonts,amssymb,amscd,amsthm}
\usepackage{epsfig}
\makeatletter \voffset=-15mm \textwidth=17cm \textheight=24cm
\newtheorem{theo}{Theorem}[section]

\newtheorem{prop}{Proposition}[section]

\makeatletter \@addtoreset{equation}{section} \makeatother

\newcommand{\mR}{\mathbb{R}}
\newcommand{\mT}{\mathbb{T}}
\newcommand{\mZ}{\mathbb{Z}}

\newcommand{\bdone}{\mathbf{1}}

\newcommand{\bdJ}{\mathbf{J}}

\newcommand{\bdc}{\mathbf{c}}

\newcommand{\calA}{\mathcal{A}}

\newcommand{\calS}{\mathcal{S}}
\newcommand{\calT}{\mathcal{T}}

\newcommand{\prob}{\mathsf P}

\newcommand{\ph}{\varphi}
\newcommand{\thet}{\vartheta}

\newcommand{\p}{\operatorname{\mbox{pr}}}

\graphicspath{{shar_pics/}}
\begin{document}

\title{On the support of a body by a surface with random roughness}
\author{D.Treschev}

\maketitle

\begin{abstract}
Suppose an interval is put on a horizontal line with random roughness. With probability one it is supported at two points, one from the left, and another from the right from its center.
We compute probability distribution of support points provided the roughness is fine grained. We also solve an analogous problem where  a circle is put on a rough plane. Some applications in static are given.
\end{abstract}

\section{Introduction}

\subsection{Motivations}

The Amonton-Coulomb law of friction (dry friction) says that if the motion of a body is a translation along a fixed plane, the friction force is up to a constant multiplier (the dry friction coefficient) equals total normal load. If the body slides along a plane with nonzero angular velocity, to obtain total friction force and total friction momentum, one has to integrate infinitesimal friction forces over the contact spot. This makes the problem of sliding of a body along a plane in the presence of dry friction non-trivial.

There is a series of publications where dynamical problems of this kind are studied: \cite{Con,Zh1,ISC,McM,Far,27_n}. A key role in these models is played by the hypothesis on the distribution of the normal load on the contact spot. All such hypotheses are essentially phenomenological although some quasistatic argument is usually attached. The uniform distribution \cite{ISC,Far,27_n} or rotational symmetric ones (for cylindrical bodies with rotational symmetric base) \cite{Kir1,Kir2} are compatible with dynamics only for bodies of infinitesimal height. Dynamically compatible deformations of the above distributions are considered in \cite{Iva}, see also \cite{GST,BS,EIT}, where qualitative analysis of the motion is presented.

Very careful experiments \cite{27_n}, where a plastic disk slides along nylon, stretched over the surface of a flat table, essentially confirm (even quantitatively) theoretical predictions. Other experiments, where a rigid disk slides along a rigid surface \cite{Kir4,GST} produce much more noisy data which correspond to the the above theoretical works only qualitatively. We believe that the main reason for such noisy and unstable data is that when both the disk and the support surface are sufficiently rigid, it is hard to expect that their surfaces are perfectly flat: very small deviations from ideal flatness can change unpredictably the distribution of the normal load and break any deterministic hypothesis on the distribution of a load over the contact spot. In this case one should use some probabilistic assumptions. For example, it is possible to consider a (perfectly) flat body on a rough surface with random roughness.

Instead of a disk on a plane in this paper we consider two simpler problems: an interval on a rough line and a circle on a rough plane. We also consider some static problems which appear in this context.

\subsection{An interval on a line}

\begin{figure}
\begin{center}
\includegraphics[scale = .7]{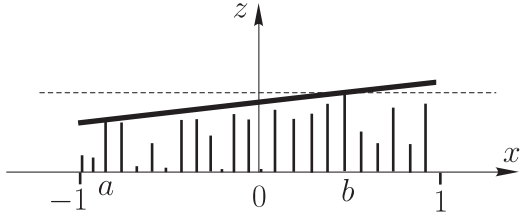}
\caption{Random comb}
\label{fig:comb}
\end{center}
\end{figure}

Consider the points
$$
  w_j = (x_j,0), \quad
  x_j = -1 + 2/N, \qquad
  j=1,\ldots,N
$$
on the horizontal interval
$$
  I = \{(x,z)\in\mR^2 : x\in [-1,1], z=0\} .
$$
Each point $w_j$ is supposed to be the lower end of a vertical interval whose length $\xi_j$ is uniformly distributed on $[0,1]$.  We call any such vertical interval a tooth and the whole set of these intervals a random comb, see Fig. \ref{fig:comb}

An interval $J$, lying on this random comb and projecting exactly on $I$, with probability 1 is supported by two teeth with horizontal coordinates
$$
  a_1 = x_{j_1}\in I_- = [-1,0], \quad
  a_2 = x_{j_1}\in I_+ = [0,1], \qquad
  1\le j_1,j_2 \le N.
$$
We say that in this case the event $S^a$ takes place.

\begin{theo}
\label{theo:int}
In the limit $N\to\infty$ the density $p : I_-\times I_+ \to \mR_+$ of probability distribution of the random event $S^a$ is
\begin{equation}
\label{p}
    p(a)
  = (a_2-a_1) \Big( \frac{4}{3(1 + a_2)^3}
                   + \frac{4}{3(1 - a_1)^3}
                    + \frac16 \Big).
\end{equation}
\end{theo}

\begin{figure}
\begin{center}
\includegraphics[scale = .6]{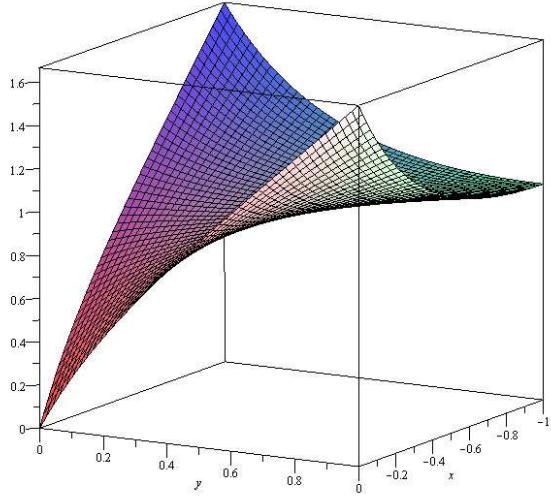}
\caption{Graph of the function p}
\label{fig:p_density}
\end{center}
\end{figure}

Graph of the function $p$ is presented in Fig. \ref{fig:p_density}. We see that $p(0,0)=0$ and $p$ attains global maximum at the points $(-1,0)$ and (0,1).

To get ``mechanical'' interpretations, suppose that the heavy interval $J$ is drawn along a rough line. Where it will be scratched more: near ends or in the middle?

Density of probability distribution for the right support point $a_2$ is as follows:
$$
    p_2
  = \int_{-1}^0 p(a)\, da_1
  = \frac{4}{3(1+a_2)^2} - \frac{2}{3(1+a_2)^3}
    + \frac{2a_2}{3} + \frac14.
$$
We have: $p_2(1)/p_2(0) = 14/11$. Therefore endpoints of $J$ are support points 14/11 times more frequently than points near the center of $J$. However we should take into account that the rate of scratching depends also on the normal load. Hence we have to perform another calculation.

Suppose that the rate of scratching is proportional to the normal load. If $J$ is supported at the points $a_1 < 0 < a_2$, the left and right tooth carries the weight
$$
  l_1(a) = \frac{P a_2}{a_2 - a_1} \quad
  \mbox{and}\quad
  l_2(a) = \frac{P a_1}{a_1 - a_2}
$$
respectively, where $P$ is the weight of $J$. Therefore the rate of scratching at the left support point is proportional to
$$
    \mbox{\rm scr}_1(a_1)
  = \int_0^1 \frac{a_2}{a_2-a_1} p(a)\, da_2
  = \frac{2}{3(1-a_1)^3} + \frac12.
$$
Analogously
$\mbox{\rm scr}_2(a_2) = \frac{2}{3(1+a_2)^3} + \frac12$.
Since $\frac{\mbox{\rm scr}_2(0)}{\mbox{\rm scr}_2(1)}
       = \frac{10}{7}$, we see that the middle point will be
scratched stronger than the end.
\medskip

\begin{figure}
\begin{center}
\includegraphics[scale = .3]{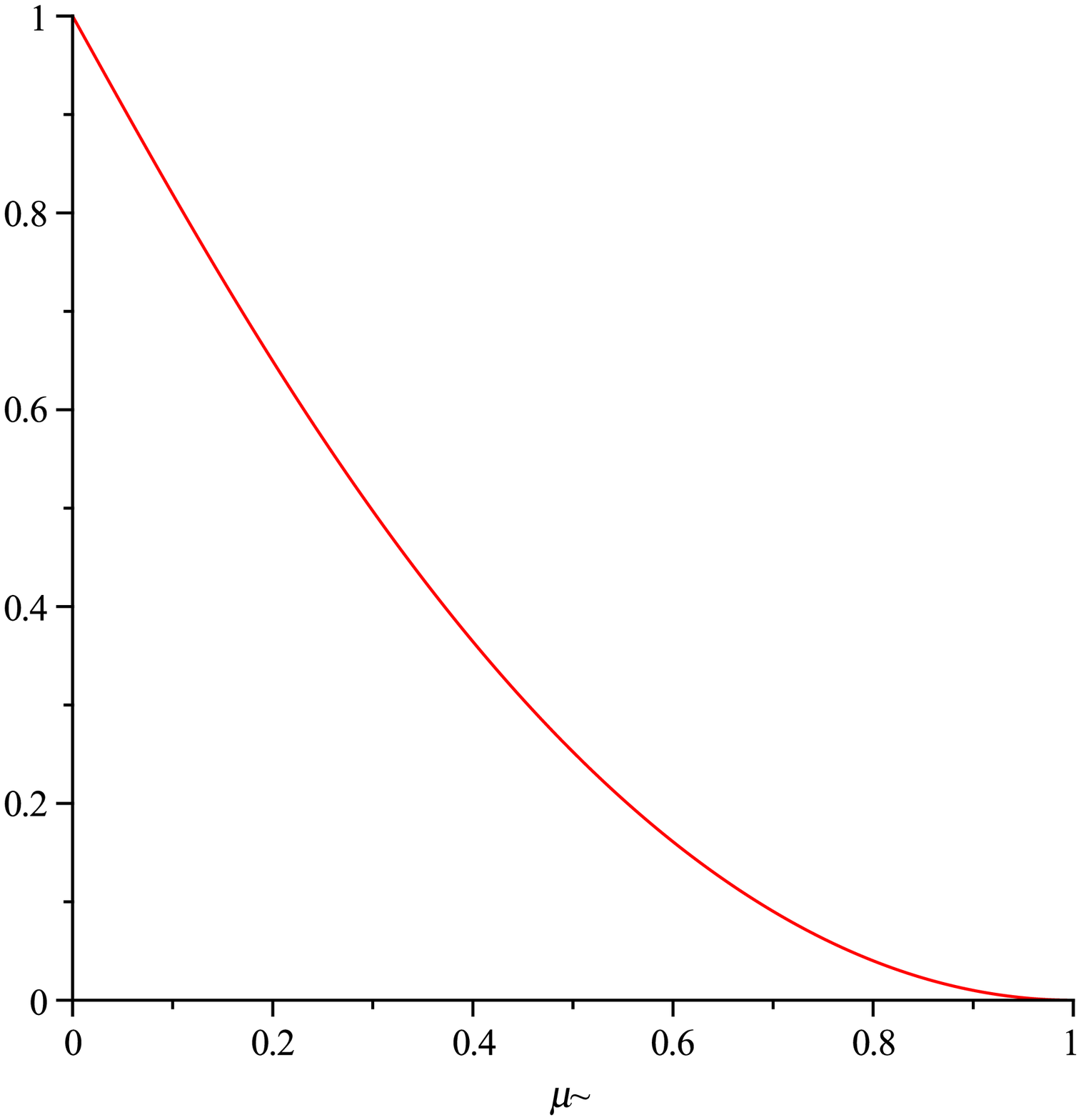}
\caption{Graph of the function $p_*$}
\label{fig:*}
\end{center}
\end{figure}

Another application of (\ref{p}) is as follows. Suppose that $J$ is a heavy beam of mass $M$ lying on an uneven surface. A man of mass $m$ walks along the beam. At some moment it may happen that under the weight of the man the beam will leave its initial equilibrium, starting to rotate on one of the support points. We compute the probability
$$
  p_* = p_*(\mu), \qquad
  \mu = \frac{m}{m+M} \in [0,1]
$$
of the random event that this does not happen. This event is equivalent to the following two inequalities:
$$
  -a_1 > \mu, \quad
   a_2 > \mu.
$$
Therefore
$$
    p_*
  = \int_{-1}^\mu \int_\mu^1
      (a_2 - a_1) \Big( \frac 4{3(1+a_2)^3}
                      + \frac 4{3(1-a_1)^3}
                      + \frac16 \Big)\, da_2 da_1
  = \frac{(1-\mu)^2}{6}\,
    \Big( 6+\mu -\Big(\frac{2\mu}{1+\mu}\Big)^2\Big) .
$$
In particular, if $m=M$, we have $p_* \approx 1/4$.
Graph of the function $p_*(\mu)$ is presented in Fig. \ref{fig:*}.

\subsection{Circle on a plane}

\begin{figure}
\begin{center}
\includegraphics[scale = .7]{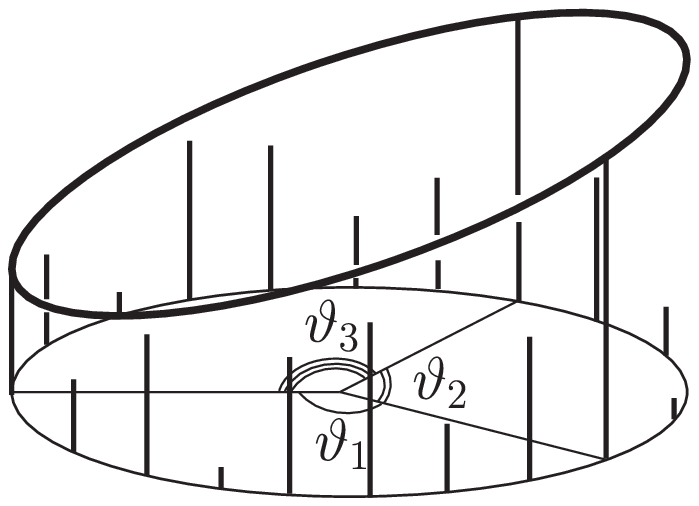}
\caption{Random circular comb}
\label{fig:can}
\end{center}
\end{figure}

Consider the points
$$
    w_j
  = w(\alpha_j)
  = (x_j,y_j,z_j)
  = (\cos\alpha_j,\sin\alpha_j,0), \quad
    \alpha_j
  = 2\pi j/N \qquad
    j = 1,\ldots,N
$$
on the horizontal circle
$$
  \bdc = \{(x,y,z)\in\mR^3 : x^2 + y^2 = 1, \; z=0 \}.
$$
Each point $w_j$ is supposed to be the lower end of a vertical interval whose length $\xi(w_j)$ is uniformly distributed on $\bdc$.  We call any such vertical interval a tooth and the whole set of these intervals a random circular comb, see Fig. \ref{fig:can}.

A thin hoop $J$, lying on this random comb, with probability 1 is supported by three teeth
$$
        w_{n_i}
      = (\cos\alpha_{n_i},\sin\alpha_{n_i},0), \qquad
  1\le n_i \le N, \quad
  i = 1,2,3.
$$
We say that in this case the event $S^\ph$ takes place
$S^\ph$ takes place, where
\begin{equation}
\label{ph=n}
  \ph = (\ph_1,\ph_2,\ph_3), \quad
  \ph_i = 2\pi n_i / N \bmod 2\pi, \qquad
  i = 1,2,3.
\end{equation}
We are interested in probability distribution of the random event $S^\ph$.

We assume that orientation of the triangle
$\ph = (\ph_1,\ph_2,\ph_3)$, $\ph_i = \alpha_{n_i}$
is positive i.e.,
$$
  \ph_{i+1} - \ph_i = \thet_{i-1}\bmod 2\pi, \quad
  \mbox{for some real }
  \thet_1,\thet_2,\thet_3 > 0, \quad
  \thet_1 + \thet_2 + \thet_3 = 2\pi,
$$
where it is convenient to assume the subscript $i$ to lie in the cyclic group $\mZ_3$.

The mass center of $J$ should lie inside the triangle with vertices $w_{n_1},w_{n_2},w_{n_3}$ (otherwise $J$ can not be in equilibrium on the teeth $n_1,n_2,n_3$). This condition is equivalent to the inequalities
$0 < \thet_i < \pi$. Moreover, the events $S^{(\ph_1,\ph_2,\ph_3)}$, $S^{(\ph_2,\ph_3,\ph_1)}$, and $S^{(\ph_3,\ph_1,\ph_2)}$ are the same. Therefore
$$
     \ph \in \widehat\calS = \calS/\mZ_3, \qquad
     \calS
  =  \{\ph\in\mT^3 : 0 < \thet_i(\ph) < \pi, \; i= 1,2,3 \},
$$
where $\mZ_3$ acts on $\mT^3$ by cyclic permutations:
$$
  (\ph_1,\ph_2,\ph_3) \mapsto (\ph_2,\ph_3,\ph_1)
                     \mapsto (\ph_3,\ph_1,\ph_2).
$$

In the limit $N\to\infty$ distribution of the random event $S^\ph$ has density $p_\calS : \widehat\calS\to\mR_+$. This density is invariant with respect to the action $R_\alpha$ of the circle $\mT$:
\begin{equation}
\label{action}
          \widehat\calS\ni\hat\ph
  \mapsto R_\alpha(\hat\ph)
      =   \hat\ph + \alpha\, \bdone, \qquad
          \bdone
      =   (1,1,1)^T\in\mR^3, \quad
          \alpha\in\mT  .
\end{equation}
Therefore it is natural to consider this distribution on the quotient
\begin{equation}
\label{triangle}
        \widehat\calT
      = \widehat\calS / \mT
      = \calT / \mZ_3, \qquad
        \calT
      = \{\thet = (\thet_1,\thet_2,\thet_3) :
         0 < \thet_i < \pi, \;
         \thet_1 + \thet_2 + \thet_3 = 2\pi \} .
\end{equation}
More precisely, let $\p : \widehat\calS\to\widehat\calT$ be the natural projection. Then there exists a function
$p_\calT : \widehat\calT\to\mR_+$ such that
$p_\calT\circ\p = p_\calS$. The space $\widehat\calT$ should be considered with the measure $\mu_\calT$:
\begin{equation}
\label{mu_T}
    d\mu_\calT
  = \frac13 \big| d\thet_3\wedge d\thet_2
                 + d\thet_1\wedge d\thet_3
                  + d\thet_2\wedge d\thet_1
            \big| .
\end{equation}
Then $d\hat\ph = d\hat\ph_1 d\hat\ph_2 d\hat\ph_3$ is the pull-back of $d\mu_\calT$: $\,\p_* (d\mu_\calT) = d\hat\ph$.

\begin{theo}
\label{theo:circ}
The density $p_\calT$ satisfies the equation
$$
    p_\calT(\hat\thet)
  = 2\pi \sin\frac{\thet_1}{2}
          \sin\frac{\thet_2}{2}
           \sin\frac{\thet_3}{2}
     \Big( \frac1{\pi^2} + \sum_{i=1}^3 f(\hat\thet_i) \Big),\quad
     f(\xi)
  =  \int_0^{\xi/2} \!
     \frac{(\xi - 2\ph)\sin\ph}
          {((\pi-\ph)\cos\ph + \sin\ph)^3}
            \, d\ph .
$$
\end{theo}

\begin{figure}
\begin{center}
\includegraphics[scale = .7]{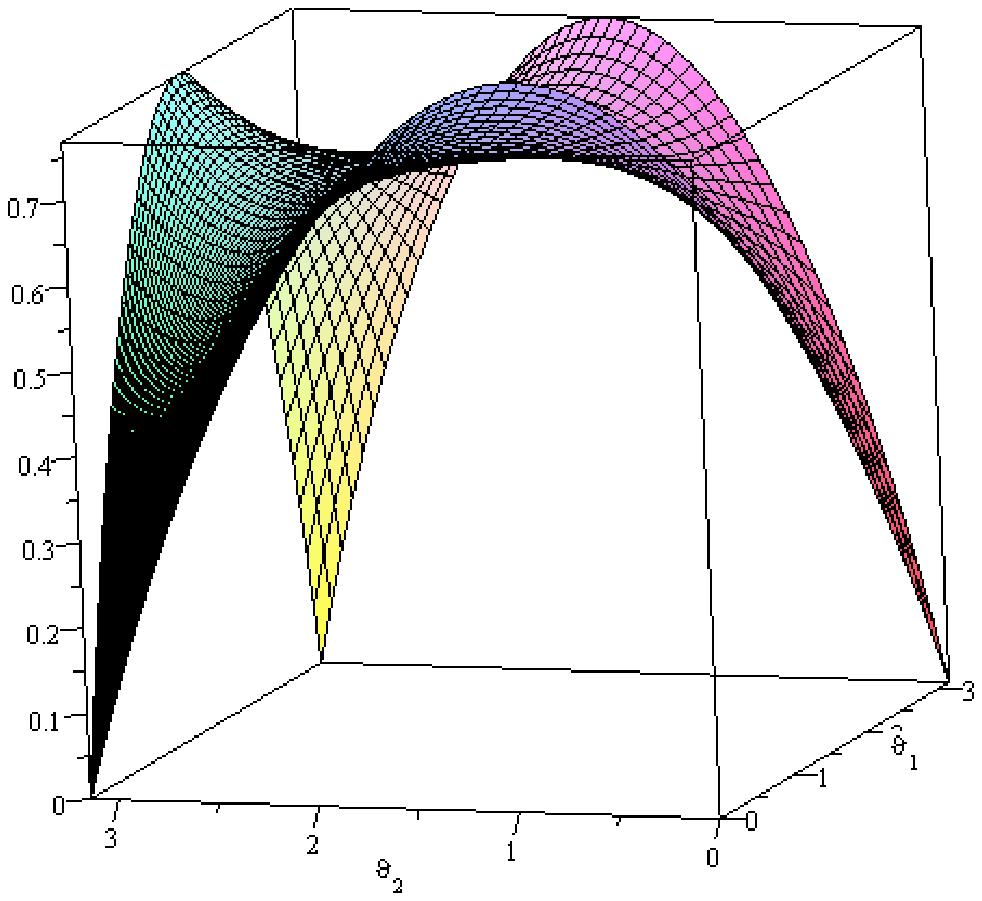}
\caption{Graph of the function $p_\calT$.}
\label{fig:comb_circ}
\end{center}
\end{figure}

Graph of the function $p_\calT$ is presented in Fig. \ref{fig:comb_circ}. Here we take $\thet_1,\thet_2$ as coordinates on $\widehat\calT$. Hence $\widehat\calT$ can be regarded as the triangle
$$
  \{(\thet_1,\thet_2) : 0<\thet_1<\pi,\,
       0<\thet_2<\pi,\, \thet_1 + \thet_2 < \pi\}
$$
with identification
\begin{equation}
\label{equivalence}
       (\thet_1,\thet_2)
  \sim (\thet_2,2\pi-\thet_1-\thet_2)
  \sim (2\pi-\thet_1-\thet_2,\thet_1).
\end{equation}

We see that $p=0$ if one of the angles $\thet_i$ vanishes. Maximal value of $p$ is attained at points $\thet$ such that for some $i\in\mZ_3$ $\; \thet_i = \pi$ and $\thet_{i\pm 1} = \pi/2$.

As an illustration consider a man of mass $m$ going around the hoop of mass $M$. Let
$$
  p_* = p_*(\mu), \qquad
  \cos\alpha = \mu = \frac{m}{m+M},\quad
  0\le\alpha\le\pi/2
$$
be the probability of the random event that the hoop stands motionless during all the walk.

This event is equivalent to the 3 inequalities
$$
  0 < \thet_i < 2\alpha, \qquad
  i = 1,2,3.
$$
Hence
$$
  p_* = \int_{D(\alpha)} p_{\calT}(\thet)\, d\thet_1 d\thet_2,
  \qquad
    D(\alpha)
  = \big\{(\thet_1,\thet_2) : \thet_1 < 2\alpha,\,
                          \thet_1 < 2\alpha,\,
                          2\pi - 2\alpha < \thet_1 + \thet_2
    \big\} / \sim ,
$$
where $\sim$ is the equivalence relation (\ref{equivalence}).

\begin{figure}
\begin{center}
\includegraphics[scale = .3]{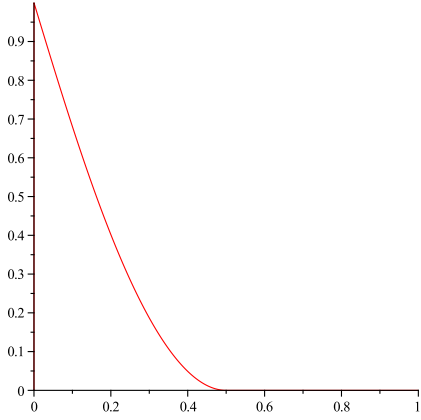}
\caption{Graph of the function $\mu\mapsto p_*(\mu)$}
\label{fig:man_on_circ}
\end{center}
\end{figure}

Since $D(\alpha)$ is empty for $\alpha < \pi/3$, we only have to consider the case $\pi/3 < \alpha < \pi/2$. Graph of the function
$\mu\mapsto p_*(\mu)$,
$$
    p_*
  = \frac{2\pi}3
    \int_{2\pi-4\alpha}^{2\alpha} d\thet_1
     \int_{2\pi-2\alpha-\thet_1}^{2\alpha}
      \sin\frac{\thet_1}{2}
       \sin\frac{\thet_2}{2}
        \sin\frac{\thet_1+\thet_2}{2}
        \Big( \frac1{\pi^2}
              + f(\thet_1) + f(\thet_2)
              + f(2\pi - \thet_1 - \thet_2) \Big)\, d\thet_2
$$
is presented in Fig. \ref{fig:man_on_circ}. In particular,
$p_* = 1/2$ for $\mu\approx 1/6$.

\section{Proof of Theorem \protect\ref{theo:int}}

Let $\Omega$ be the configuration space of the random comb:
$\Omega = [0,1])^N$. We consider large integer $L$, and put
$$
  n = (n_1,n_2), \quad
  1\le n_1 \le N/2 < n_2 \le N.
$$
Then we define two random events $\nu$ and $Q_n$, where by definition
\begin{itemize}
\item $\nu = n$ iff $J$ is supported by the teeth $n_1$ and
      $n_2$,
\item $Q_n = (K_1,K_2)$ iff length of the tooth with number
      $n_i$ equals
      $$
            \xi_{n_i}
        \in \big(1 - (K_i-1) / (NL),1 - K_i / (NL)\big),
            \qquad  i = 1,2.
      $$
\end{itemize}

For any $K\in \{1,\ldots,NL\}^2$ we have:
$$
  \prob\{Q_n = K\} = (NL)^{-2}.
$$
Therefore by the formula of total probability
$$
    \prob\{\nu = n\}
  = \sum_K \prob\{\nu = n | Q_n = K\}\, \prob\{Q_n=K\}
  = \sum_K \frac{\prob\{\nu = n | Q_n = K\}}
                {(NL)^2}.
$$
In the limit $L\to\infty$ we obtain:
\begin{equation}
\label{limit}
    \prob\{\nu=n\}
  = \frac{1}{N^2}
    \int_0^N\!\!\!\int_0^N
      \prob\{\nu=n | \xi_n = \bdone - A/N\} \, dA, \qquad
      \xi_n
  = (\xi_{n_1},\xi_{n_2}), \quad
    \bdone = (1,1)^T .
\end{equation}

In the limit $N\to\infty$ we obtain densities of probability distributions
\begin{eqnarray*}
&   p : I_-\times I_+ \to \mR_+, \quad
    p_{\nu|Q}:I_-\times I_+\times [0,N]^2 \to\mR_+, & \\
&\displaystyle
    p(a)
  = \lim_{N\to\infty} \frac{N^2}{4}\,\prob\{\nu=n\}, \quad
    p_{\nu|Q}(a,A)
  = \lim_{N\to\infty} \prob\{\nu = n, \xi_n = 1-A/N\},  & \\
&   a = (a_1,a_2) = \big(-1 + 2n_1/N,-1 + 2n_2/N\big).  &
\end{eqnarray*}

Equation (\ref{limit}) implies
\begin{equation}
\label{p(a)}
    p(a)
  = \frac14 \int_0^N\!\!\int_0^N p_{\nu|Q}(a,A)\, dA .
\end{equation}

Now we turn to computation of $p_{\nu|Q}$. The interval
$J = J(a,A)$ is determined by the equation
$$
  z = 1 - \frac{A_1a_2-A_2a_1}{N(a_2-a_1)}
        - \frac{A_2-A_1}{N(a_2-a_1)}x, \qquad
  x\in I.
$$
We have to consider two cases.
\medskip

(1) The interval $J$ does not intersect the line segment $I_+$ joining the points $(-1,1)$ and $(1,1)$. This happens provided
$$
     |A_1-A_2| < A_1a_2 - A_2a_1.
$$

(2) $J\cap I_+ = (x_*,1)$. In this case
$|A_1-A_2| \ge A_1a_2-A_2a_1$ and
$$
  x_* = \frac{A_1a_2-A_2a_1}{A_1-A_2}.
$$

In case (1) probability for the point $w_j=(x_j,0)$ to have the tooth (entirely) under $J$ is
\begin{equation}
\label{under}
    z_j
  = 1 - \frac{A_1a_2-A_2a_1}{N(a_2-a_1)}
      - \frac{A_2-A_1}{N(a_2-a_1)} x_j
  < 1.
\end{equation}

Therefore probability for the whole comb to be under $J$ is
\begin{eqnarray*}
&\displaystyle
    p_{\nu|Q}^{(1)}
  = \prod_{j=1}^N
        \Big( 1 - \frac{A_1a_2-A_2a_1}{N(a_2-a_1)}
                - \frac{A_2-A_1}{N(a_2-a_1)}
                     \Big(-1 + \frac{2j}{N}\Big)
        \Big)
  = e^{F_1}, & \\
&\displaystyle
    F_1
  = \sum_{j=1}^N
        \log \Big( 1 - \frac{A_1a_2-A_2a_1}{N(a_2-a_1)}
                - \frac{A_2-A_1}{N(a_2-a_1)}
                     \Big(-1 + \frac{2j}{N}\Big)
             \Big). &
\end{eqnarray*}
In the limit $N\to\infty$ we have:
$F_1 = - \frac{A_1a_2-A_2a_1}{a_2-a_1}$. Hence
$$
  p_{\nu|Q}^{(1)} = e^{- \frac{A_1a_2-A_2a_1}{a_2-a_1}} .
$$

Consider case (2). For definiteness we assume that $A_1>A_2$ i.e.,
$x_*>0$. Then probability for any point $w_j$ to have the tooth under $J$ is determined by (\ref{under}) if $x_j\in [-1,x_*]$ and equals 1 if $x_j\in [x_*,1]$.

Probability for the whole comb to be under $J$ is
$p^{(2)}_{\nu|Q} = e^{F_2}$,
\begin{eqnarray*}
     F_2
 &=& \sum_{j\ge 1, -1+2j/N\le x_*}
        \log \Big( 1 - \frac{A_1a_2-A_2a_1}{N(a_2-a_1)}
                - \frac{A_2-A_1}{N(a_2-a_1)}
                     \Big(-1 + \frac{2j}{N}\Big)
             \Big) \\
 &=& -\frac12 \int_{-1}^{x_*}
           \Big( \frac{A_1a_2-A_2a_1}{a_2-a_1}
                + \frac{A_2-A_1}{a_2-a_1}x
           \Big)\, dx + O(1/N) .
\end{eqnarray*}
For $N\to\infty$ we obtain:
$$
    p_{\nu|Q}^{(2)}
  = e^{- \frac{(A_1(a_2+1)-A_2(a_1+1))^2}{4(a_2-a_1)(A_1-A_2)}}
  \qquad \mbox{if $A_1>A_2$}.
$$
The case $A_1<A_2$ can be obtained from this one by the exchange
$A_1\leftrightarrow A_2$, $a_1\leftrightarrow -a_2$. Therefore
$$
    p_{\nu|Q}^{(2)}
  = e^{- \frac{(A_2(-a_1+1)-A_1(-a_2+1))^2}{4(a_2-a_1)(A_2-A_1)}}
  \qquad \mbox{if $A_2>A_1$}.
$$

Considering in (\ref{p(a)}) the limit $N\to\infty$ we see that
\begin{eqnarray}
\label{pQQQ}
&\displaystyle
    p(a)|_{N\to\infty} = Q_1 + Q_2^+ + Q_2^-,  \qquad
       Q_1 = \frac14 \int_{D_1} p_{\nu|Q}^{(1)}\, dA, \quad
   Q_2^\pm = \frac14 \int_{D_2^\pm} p_{\nu|Q}^{(2)}\, dA ,& \\
\nonumber
&    D_1
   = \{ A \in \mR_+^2 : |A_1-A_2| < A_1a_2-A_2a_1 \},  & \\
\nonumber
&    D_2^\pm
   = \{ A \in \mR_+^2 : |A_1-A_2| \ge A_1a_2-A_2a_1, \;
                           \pm(A_1-A_2) > 0 \}. &
\end{eqnarray}

Change of the variables
$$
  a_2A_1 - a_1A_2 = r, \quad - A_1 + A_2 = q
$$
transforms the integrals as follows:
\begin{eqnarray}
\label{Q1}
     Q_1
 &=& \frac14 \int_{|q|<r} \frac{1}{a_2-a_1}
                    e^{- \frac{r}{a_2-a_1}}\, drdq
  =  \frac12 (a_2-a_1), \\
\nonumber
     Q_2^+
 &=& \frac14 \int_{D^+} \frac{1}{a_2-a_1}
                  e^{\frac{r^2}{4(a_2-a_1)q}} \, drdq, \qquad
     D^+
  =  \{ -(1+a_2)q < r < -2q \}.
\end{eqnarray}
The quantity $Q_2^-(a)$ is obtained from $Q_2^+(a)$ by the exchange $a_1\leftrightarrow -a_2$.

It is convenient to compute $Q_2^+$ in the variables $u = -r^2/q$, $v = -r/q$. Direct computation gives:
\begin{equation}
\label{Q2}
    Q_2^+
  = \frac{4(a_2-a_1)}{3}
      \Big(\frac{1}{(1+a_2)^3} - \frac18\Big),\quad
    Q_2^-
  = \frac{4(a_2-a_1)}{3}
      \Big(\frac{1}{(1-a_1)^3} - \frac18\Big) .
\end{equation}
Now equation (\ref{p}) follows from (\ref{pQQQ}), (\ref{Q1}), and
(\ref{Q2}). \qed

\section{Proof of Theorem \protect\ref{theo:circ}}

Let $\Omega$ be the configuration space of the circular random comb: $\Omega = [0,1]^N$. The teeth that support $J$ are determined by
equation (\ref{ph=n}).

We consider large integer $L$, and define two random events $\nu$ and $Q_n$, where by definition
\begin{itemize}
\item $\nu = n$ iff $J$ is supported by the teeth
      $n = (n_1,n_2,n_3)$,
\item $Q_n = K = (K_1,K_2,K_3)$ iff length of the tooth with
      number $n_i$ equals
      $$
            \xi_{n_i}
        \in \big(1 - (K_i-1) / (NL),1 - K_i / (NL)\big),
            \qquad  i = 1,2,3.
      $$
\end{itemize}

For any $K\in \{1,\ldots,NL\}^3$ we have:
$\prob\{Q_n = K\} = (NL)^{-3}$.
Therefore by the formula of total probability
$$
    \prob\{\nu = n\}
  = \sum_K \prob\{\nu = n | Q_n = K\}\, \prob\{Q_n=K\}
  = \sum_K \frac{\prob\{\nu = n | Q_n = K\}}
                {(NL)^3}.
$$
Putting $\bdone = (1,1,1)^T\in\mR^3$, in the limit $L\to\infty$ we obtain:
\begin{equation}
\label{limit_o}
    \prob\{\nu=n\}
  = \frac{1}{N^3}
    \int_0^N\!\!\!\int_0^N\!\!\!\int_0^N
      \prob\{\nu=n | \xi_n = \bdone - A/N\} \, dA, \qquad
      \xi_n
  = (\xi_{n_1},\xi_{n_2},\xi_{n_3}) .
\end{equation}

In the limit $N\to\infty$ we obtain densities of probability distributions
\begin{eqnarray}
\nonumber
&   p_\calS : \widehat\calS \to \mR_+, \quad
    \tilde p_{\nu|Q} : \calS\times [0,N]^3 \to\mR_+,  & \\
\label{limlim}
&\displaystyle
    p_\calS(\hat\ph)
  = \lim_{N\to\infty}
      \Big(\frac{N}{2\pi}\Big)^3\,\prob\{\nu=n\}, \quad
    \tilde p_{\nu|Q}(\ph|A)
  = \lim_{N\to\infty} \prob\{\nu = n | \xi_n = 1-A/N\},  &
\end{eqnarray}
where $\hat\ph\in\widehat\calS$ and $\ph\in\calS$.

Equations (\ref{limit_o})--(\ref{limlim}) imply
\begin{equation}
\label{p(a)o}
    p_\calS(\hat\ph)
  = \frac1{8\pi^3} \int_{\mR^3_+}
                  \tilde p_{\nu|Q}(\ph|A)\, dA , \qquad
    \mR_+^3
  = \{A = (A_1,A_2,A_3)\in\mR^3 : A_i > 0,\; i=1,2,3\}  ,
\end{equation}
where $\hat\ph$ is the image of $\ph$ under the natural map
$\calS\to\widehat\calS$.

Both densities $p_\calS$ and $\tilde p_{\nu|Q}$ are invariant with respect to the action $R_\alpha$ of the group $\mT$, see (\ref{action}). Hence we obtain the densities $p_\calT,p_{\nu|Q}$ on $\widehat\calT = \widehat\calS/\mT$ and
$\calT\times [0,N]^3$ respectively: $\calT = \calS / \mT$,
$$
  p_\calT(\hat\thet) = 2\pi p_\calS(\hat\ph), \quad
  p_{\nu|Q}(\thet|A) = \tilde p_{\nu|Q}(\ph|A),
$$
where measures on $\hat\calT$ and $\calT$ are determined by (\ref{mu_T}). Then (\ref{p(a)o}) implies
\begin{equation}
\label{p(thet)}
    p_\calT(\hat\thet)
  = \frac1{4\pi^2}
      \int_{\mR_+^3}
         p_{\nu|Q}(\thet|A)\, dA.
\end{equation}
Now we turn to computation of $p_{\nu|Q}$.

The plane passing through $J = J(a(\ph),A)$ is determined by the equation
\begin{eqnarray}
\nonumber
&\displaystyle
  z = 1 - \frac{\sigma_0}{N}
        - \frac{\sigma_x}{N}x
        - \frac{\sigma_y}{N}y , & \\
\label{sigma0}
&\displaystyle
    \sigma_0
  = \frac{1}{\Delta}
     \left|\begin{array}{ccc}
       \cos\ph_1 & \sin\ph_1 & A_1 \\
       \cos\ph_2 & \sin\ph_2 & A_2 \\
       \cos\ph_3 & \sin\ph_3 & A_3
     \end{array}\right|
  = \frac{A_1\sin\thet_1 + A_2\sin\thet_2 + A_3\sin\thet_3}
         {\Delta} > 0, & \\
\label{s+s+s}
&\displaystyle
    \Delta
  = \left|\begin{array}{ccc}
       \cos\ph_1 & \sin\ph_1 & 1 \\
       \cos\ph_2 & \sin\ph_2 & 1 \\
       \cos\ph_3 & \sin\ph_3 & 1
     \end{array}\right|
  = \sin\thet_1 + \sin\thet_2 + \sin\thet_3, & \\
\label{sigmaxy}
&\displaystyle
    \sigma_x
  = \frac{1}{\Delta}
     \left|\begin{array}{ccc}
       A_1 & \sin\ph_1 & 1 \\
       A_2 & \sin\ph_2 & 1 \\
       A_3 & \sin\ph_3 & 1
     \end{array}\right|, \quad
    \sigma_y
  = \frac{1}{\Delta}
    \left|\begin{array}{ccc}
       \cos\ph_1 & A_1 & 1 \\
       \cos\ph_2 & A_2 & 1 \\
       \cos\ph_3 & A_3 & 1
     \end{array}\right| .
\end{eqnarray}

We consider two cases.
\medskip

(1) The disk $J$ does not intersect the disk $I_+$, obtained as a shift of the disk $I$ by the vector $(0,0,1)$. This happens provided $\sigma_x\cos\ph + \sigma_y\sin\ph + \sigma_0 > 0$ for all real $\ph$ i.e.,
$$
  \sigma_x^2 + \sigma_y^2 < \sigma_0^2.
$$

(2) $J\cap I_+ \ne \emptyset$. In this case $J$ is below $I_+$ over the domain
$$
   D_\sigma
 = \{x^2 + y^2 \le 1, \;
      \sigma_x x + \sigma_y y + \sigma_0 \ge 0 \}.
$$

In case (1) probability for the point $w_j = w(\alpha_j)$ to have a tooth (entirely) under $J$ is
$$
        1 - \frac{\sigma_0}{N}
          - \frac{\sigma_x}{N} \cos\alpha_j
          - \frac{\sigma_y}{N} \sin\alpha_j
     \le 1.
$$
Therefore probability for the whole comb to be under $J$ equals
\begin{eqnarray*}
&   p_A^{(1)}
  = \prod_{j\ne n_1,n_2,n_3}
         \Big( 1 - \frac{\sigma_0}{N}
                  - \frac{\sigma_x}{N} \cos\alpha_j
                   - \frac{\sigma_y}{N} \sin\alpha_j \Big)
  = e^{F_1}, \qquad
    \alpha_j = \frac{2\pi j}{N}, \quad
    j = 1,\ldots,N, & \\
&   F_1
  = \sum_{j\ne n_1,n_2,n_3}
     \log\Big( 1 - \frac{\sigma_0}{N}
                  - \frac{\sigma_x}{N} \cos\alpha_j
                   - \frac{\sigma_y}{N} \sin\alpha_j \Big)
  = - \sigma_0 + O(1/N). &
\end{eqnarray*}
For $N\to\infty$ we obtain:
$$
  p_A^{(1)} = e^{-\sigma_0} .
$$

In case (2) the tooth is under $J$ with probability
\begin{eqnarray*}
&\displaystyle
   1 - \frac{\sigma_0
             + \sigma_x\cos\alpha_j
              + \sigma_y\sin\alpha_j}{N} \quad
   \mbox{if } \alpha_j \in B^+, \quad
   \mbox{and $1$ if } \alpha_j \in B^-, & \\[1mm]
&   B^\pm
  = \big\{\alpha\in\mT :
           \pm (\sigma_x\cos\alpha
                + \sigma_y\sin\alpha + \sigma_0) \ge 0 \big\}. &
\end{eqnarray*}
Therefore probability for the whole comb to be under $J$ equals
\begin{eqnarray*}
&   p_A^{(2)}
  = \prod_{\alpha_j\in B^+,\, j\ne n_1,n_2,n_3}
      \Big(1 - \frac{\sigma_0
             + \sigma_x\cos\alpha_j
              + \sigma_y\sin\alpha_j}{N}
      \Big)
  = e^{F_2}, \qquad
    \alpha_j = \frac{2\pi j}{N}, \quad
    j = 1,\ldots,N, & \\
&   F_2
  = \sum_{\alpha_j\in B^+,\, j\ne n_1,n_2,n_3}
     \log \Big(1 - \frac{\sigma_0
                  + \sigma_x\cos\alpha_j
                   + \sigma_y\sin\alpha_j}{N}
          \Big)
  = \sigma_0\calA + O(1/N), & \\
&\displaystyle
    \calA
  = \frac{1}{2\pi\sigma_0}
     \int_{B^+} (\sigma_x\cos\ph
                      + \sigma_y\sin\ph + \sigma_0)\, d\ph. &
\end{eqnarray*}

\begin{prop}
\label{prop:calA}
$\calA = \frac1\pi (\ph_\sigma - \tan\ph_\sigma)$, where
$\ph_\sigma
 = \arccos\big(-\sigma_0 / \sqrt{\sigma_x^2 + \sigma_y^2}\big)$.
\end{prop}

In the limit $N\to\infty$ we obtain the probability
$$
  p_A^{(2)} = e^{-\sigma_0\calA} .
$$

By (\ref{p(thet)}) we have the equation
\begin{eqnarray}
\label{ppp}
&\displaystyle
  p = p_1 + p_2, \qquad
  p_1 = \frac1{4\pi^2}\int_{B_1} p_A^{(1)}\, dA, \quad
  p_2 = \frac1{4\pi^2}\int_{B_2} p_A^{(2)}\, dA, & \\
\nonumber
&  B_1
 = \big\{ A\in\mR^3_+ : \sigma_x^2 + \sigma_y^2 < \sigma_0^2
   \big\},\quad
   B_2
 = \big\{ A\in\mR^3_+ : \sigma_x^2 + \sigma_y^2 \ge \sigma_0^2
   \big\} . &
\end{eqnarray}
Computation of the integrals $p_1,p_2$ requires some preliminary work. First, we introduce new coordinates
$$
  \tau_i = (1 - \cos\thet_i) A_i / \Delta, \qquad i=1,2,3
$$
and put
$$
  \bdone = (1,1,1)^T \in\mR^3, \quad
  c_i = \cot \frac{\thet_i}{2}, \qquad i=1,2,3 .
$$

\begin{prop}
\label{prop:iden}
For any $\thet\in\calT$
\begin{equation}
\label{iden}
    c_1c_2 + c_2c_3 + c_3c_1
  = 1, \quad
    \frac{\sum\sin\thet_i}{\prod\sin(\thet_i/2)}
  = 4, \quad
    \frac{1}{\prod\sin(\thet_i/2)}
  = \langle\bdone,c\rangle - c_1c_2c_3.
\end{equation}
\end{prop}

Combining (\ref{sigma0})--(\ref{sigmaxy}) and (\ref{iden}), we have:
\begin{eqnarray}
\label{Delta}
     \Delta
 &=& \frac{4}{\langle c,\bdone\rangle - c_1c_2c_3} , \quad
     \sigma_0
\;=\; \langle c,\tau\rangle,   \\
\label{J}
     \sigma_x^2 + \sigma_y^2 - \sigma_0^2
 &=& -\langle\tau, J\tau\rangle, \qquad
     \bdJ
  =  \left(\begin{array}{rrr}
            -1 &  1 &  1 \\
             1 & -1 &  1 \\
             1 &  1 & -1
           \end{array}\right)
  = \bdone\otimes\bdone - 2 .
\end{eqnarray}

Integrals (\ref{ppp}) in the new coordinates take the form
\begin{eqnarray}
\nonumber
&\displaystyle\!\!\!\!\!\!
     p_i
  =  \frac{\Delta^3\hat p_i}
          {32\pi^2\prod\sin^2\frac{\thet_i}2}
  =  \frac{2\hat p_i}
          {\pi^2(\langle c,\bdone\rangle - c_1c_2c_3)}
  =  \frac{2}{\pi^2} \sin\frac{\thet_1}{2}
                      \sin\frac{\thet_2}{2}
                       \sin\frac{\thet_3}{2},     &  \\
\label{hatp_i}
&\displaystyle\!\!\!\!\!\!
     \hat p_1
  =  \int_{C_1} e^{-\langle c,\tau\rangle} \, d\tau,
     \quad
     \hat p_2
  =  \int_{C_2}
        e^{-\langle c,\tau\rangle\,\calA} \, d\tau,
        &\\[1mm]
\label{C2}
&    C_1
  =  \{ \tau\in\mR^3 : \langle c,\tau \rangle > 0, \;
              \langle\tau,\bdJ\tau\rangle > 0 \} , \quad
     C_2
  =  \{ \tau\in\mR_+^3 : \langle\tau,\bdJ\tau\rangle < 0 \} . &
\end{eqnarray}

\subsection{Convenient variables}

To compute integrals (\ref{hatp_i}), it is convenient to introduce new variables. We put
$$
  w = 1 - \frac{\langle \tau,\bdJ\tau\rangle}
               {\langle c,\tau\rangle^2}, \qquad
        \lambda
    =   \frac{|c\times\bdone|}{\sqrt 2}
 \equiv \sqrt{\langle c,\bdone\rangle^2 - 3}.
$$
Equation (\ref{J}) implies
\begin{equation}
\label{sigsig}
  (\sigma_x^2 + \sigma_y^2) / \sigma^2 = w^2 .
\end{equation}
In the domain $C_2$ (see (\ref{C2})) we have: $w > 1$. The identity
$$
    \lambda^2 w
  = \Big(\frac{\langle \bdone,\tau\rangle}
              {\langle c,\tau\rangle}
         - \langle c,\bdone\rangle
    \Big)^2
   + \frac{\langle c\times\bdone,\tau\rangle^2}
          {\langle c,\tau\rangle^2}
$$
suggests the following change of variables:
$(\tau_1,\tau_2,\tau_3)\mapsto (u,w,\psi)$,
$$
  \langle c,\tau\rangle = u, \quad
    \frac{1}{\lambda}\,
      \Big( \frac{\langle \bdone,\tau\rangle}
                 {\langle c,\tau\rangle}
           - \langle c,\bdone\rangle  \Big)
  = \sqrt{w}\cos\psi , \quad
    \frac{\langle c\times\bdone,\tau\rangle}
         {\lambda \langle c,\tau\rangle}
  = \sqrt{w}\sin\psi.
$$
The Jacobian $\det\frac{\partial(u,w,\psi)}{\partial(\tau_1,\tau_2,\tau_3)}$
equals the product
$$
  \det\frac{\partial(u,w,\psi)}{\partial(u,v_1,v_2)}
   \det\frac{\partial(u,v_1,v_2)}{\partial(\tau_1,\tau_2,\tau_3)},
  \qquad
   v_1 = \sqrt{w}\cos\psi, \quad
   v_2 = \sqrt{w}\sin\psi.
$$
These determinants equal 2 and $2u^{-2}$ respectively. Therefore
$$
   \det\frac{\partial(u,w,\psi)}{\partial(\tau_1,\tau_2,\tau_3)}
 = \frac{4}{u^2}.
$$

Assuming $i$ to be an element of the cyclic group $\mZ_3$, we put
$$
    a_i
  = \frac{c_i\langle c,\bdone\rangle - c^2}
         {(c_{i-1}+c_{i+1})\lambda},
    \quad
    b_i
  = \frac{c_{i-1} - c_{i+1}}{(c_{i-1}+c_{i+1})\lambda}.
$$
Direct computations show that
$$
  a_i^2 + b_i^2 = 1, \quad
    a_{i+1} b_{i-1} - a_{i-1} b_{i+1}
  = \frac{2c_i}{c_i^2+1}
  = \sin\thet_i, \quad
    b_{i+1} b_{i-1} + a_{i+1} a_{i-1}
  = \frac{c_i^2-1}{c_i^2+1}
  = \cos\thet_i.
$$
Therefore for some $\psi_1,\psi_2,\psi_3\in\mT$
\begin{equation}
\label{ab}
  a_i = \sin\psi_i, \quad
  b_i = \cos\psi_i, \qquad
  \psi_{i+1} - \psi_{i-1} = \thet_i.
\end{equation}

\subsection{The integrals $\hat p_1$ and $\hat p_2$}

By using the variables $(u,v,\psi)$ in (\ref{hatp_i}), we obtain:
\begin{equation}
\label{intduint}
    \hat p_1
  = \int_0^\infty du \int_{G_1} \frac{u^2}{4} e^{-u}\, dwd\psi,
    \quad
    \hat p_2
  = \int_0^\infty du \int_{G_2}
         \frac{u^2}{4} e^{-u\calA}\, dwd\psi,
\end{equation}
where the domains $G_1,G_2$ are as follows:
$$
     G_1
  =  \big\{ (w,\psi) : 0 < w < 1 \big\} , \quad
     G_2
  = \Big\{ (w,\psi) : w > 1, \;
            \frac{\tau_i}{\langle c,\tau\rangle} > 0, \;
            i = 1,2,3
    \Big\} .
$$

\begin{prop}
\label{prop:p1p2}
\begin{equation}
\label{p1p2}
    \hat p_1 = \pi , \quad
    \hat p_2
  = \pi^3 \big( f(\thet_1) + f(\thet_2) + f(\thet_3) \big) .
\end{equation}
\end{prop}

{\it Proof of Proposition \ref{prop:p1p2}}. The first equation (\ref{p1p2}) is obvious. To prove the second one, we note that
\begin{eqnarray}
\label{hatp2}
&\displaystyle
     \hat p_2
  =  \frac12
      \int_{w>1} dw \int_{\psi\in G(w)} \calA^{-3} \, d\psi , &\\
\label{G}
& G(w) = \{\psi\in\mT : \sin(\psi+\psi_i) < 1/\sqrt{w}, \quad
             i\in\mZ_3\}. &
\end{eqnarray}
Indeed, by Proposition \ref{prop:calA} and equation (\ref{sigsig})
we have: $\calA = \calA(w)$. Therefore we can perform integration in (\ref{intduint}) in the variable $u$ which implies (\ref{hatp2}).

To check that the domain $G$ is determined by (\ref{G}), we define
\begin{eqnarray*}
      \nu
  &=& \langle \bdone,\tau\rangle
   =  u(\lambda\sqrt{w} \cos\psi + \langle c,\bdone\rangle), \\
      \beta
  &=& \langle c\times\bdone,\tau\rangle
   =  u\lambda \sqrt{w} \sin\psi.
\end{eqnarray*}
Then
\begin{eqnarray*}
    \left( \begin{array}{c}
           u\\ \nu\\ \beta
           \end{array}\right)
&=& \left( \begin{array}{ccc}
           c_1    & c_2     & c_3 \\
            1     &  1      &  1  \\
          c_2-c_3 & c_3-c_1 & c_1 - c_2
           \end{array}\right)
    \left( \begin{array}{c}
           \tau_1\\ \tau_2\\ \tau_3
           \end{array}\right) , \\
    \left( \begin{array}{c}
           \tau_1\\ \tau_2\\ \tau_3
           \end{array}\right)
&=& \frac{1}{\lambda^2}
    \left( \begin{array}{ccc}
      3c_1 - \langle c,\bdone\rangle
          & c^2 - c_1\langle c,\bdone\rangle & c_2 - c_3 \\
      3c_2 - \langle c,\bdone\rangle
          & c^2 - c_2\langle c,\bdone\rangle & c_3 - c_1 \\
      3c_3 - \langle c,\bdone\rangle
          & c^2 - c_3\langle c,\bdone\rangle & c_1 - c_2
           \end{array}\right)
    \left( \begin{array}{c}
           u\\ \nu\\ \beta
           \end{array}\right)
\end{eqnarray*}
Hence the inequalities $\tau_i>0$ take the form
$$
     3c_i - \langle c,\bdone\rangle
   + (c^2 - c_i\langle c,\bdone\rangle)
       (\lambda\sqrt{w}\cos\psi + \langle c,\bdone\rangle)
   + (c_{i+1} - c_{i-1}) \lambda\sqrt{w} \sin\psi
  > 0.
$$
After simple transformations we get:
$a_i\cos\psi + b_i\sin\psi < 1/\sqrt{w}$
which implies (\ref{G}).
\medskip

Equations (\ref{hatp2})--(\ref{G}) imply
$$
     \hat p_2
  =  \frac12 \int_{w > 1} \frac{|G(w)|}{\calA^3(w)}\, dw ,
$$
where $|G(w)|$ is the measure of the set $G(w)$.

The set $\mT\setminus\{\pi/2-\psi_1,\pi/2-\psi_2,\pi/2-\psi_3\}$
has 3 connected components: $U_1,U_2$, and $U_3$, where the interval $U_i$ has endpoints $\pi/2 - \psi_{i-1}$ and $\pi/2 - \psi_{i+1}$.
Hence
\begin{eqnarray*}
& G(w) = G_1(w) + G_2(w) + G_3(w), \qquad
  G_i(w) = G\cap U_i, & \\
& \hat p_2 = \hat p_2^{(1)} + \hat p_2^{(2)} + \hat p_2^{(3)},
  \qquad
   \hat p_2^{(i)}
 = \int_{w > 1} \frac{|G_i(w)|}{\calA^3(w)} \, dw . &
\end{eqnarray*}

By using (\ref{ab}) and (\ref{G}), we get:
$$
    |G_i(w)|
  = \thet_i - 2(\pi - \ph_\sigma(w))
$$
provided the right-hand side is non-negative.

By using the change $w = 1 / \cos^2\ph_\sigma$,
$\ph_\sigma\in (\pi - \thet_i/2,\pi)$ in the integral
$$
    \hat p_2^{(i)}
  =  \frac{\pi^3}{2} \int_{w>1}
     \frac{\thet_i - 2(\pi - \ph_\sigma(w))}
          {(\ph_\sigma(w) - \tan\ph_\sigma(w))^3} \, dw ,
$$
we obtain the equation
$$
    \hat p_2^{(i)}
  = \int_{\pi - \thet_i/2}^\pi
     \frac{\pi^3 (\thet_i - 2(\pi - \ph_\sigma))\sin\ph_\sigma}
          {(-\ph_\sigma\cos\ph_\sigma + \sin\ph_\sigma)^3}
            \, d\ph_\sigma
  = 2\pi^3 f(\thet_i) .
$$

\section{Several proofs}

{\bf Proof of Proposition \ref{prop:calA}}. We define $\sigma_*$ and $\ph_*$ by the equations
$$
  \sqrt{\sigma_x^2 + \sigma_y^2} = \sigma_*, \quad
  \cos\ph_* = \sigma_x / \sigma_*, \quad
  \sin\ph_* = \sigma_y / \sigma_*
$$
Then $\calA$ takes the form
\begin{eqnarray*}
&   \calA
  = \frac{1}{2\pi}
     \int_{\hat B_+}
      \Big(1 + \frac{\sigma_*}{\sigma_0} \cos(\ph-\ph_*)
      \Big)\,d\ph, & \\
&   \hat B_+
  = \{\ph\in\mT : \sigma_0 + \sigma_*\cos(\ph-\ph_*) \ge 0\}
  = \{\ph\in\mT : -\ph_\sigma\le\ph-\ph_*\le\ph_\sigma\}. &
\end{eqnarray*}
This implies the required assertion. \qed

{\bf Proof of Proposition \ref{prop:iden}}.
The first identity (\ref{iden}) follows from the equation
$$
     c_3
  = -\cot(\thet_1/2 + \thet_2/2)
  = \frac{1 - c_1c_2}{c_1 + c_2}.
$$
To prove the second one we note that
$$
    \sin\thet_1
  = -2\sin(\thet_1/2)\cos(\thet_2/2+\thet_3/2)
  = 2\prod\sin(\thet_i/2)
     - 2\sin(\thet_1/2)\cos(\thet_2/2)\cos(\thet_3/2).
$$
Adding to this equation two analogous ones and dividing by
$\prod\sin(\thet_i/2)$, we get:
$$
    \frac{\sum\sin\thet_i}{\prod\sin(\thet_i/2)}
  = 6 - 2(c_2c_3 + c_3c_1 + c_1c_2)
  = 4.
$$

Finally, adding up the equations
\begin{eqnarray*}
     \sin^2\frac{\thet_1}2
 &=& \sin\frac{\thet_1}2 \sin\frac{\thet_2+\thet_3}2
\;=\; \sin\frac{\thet_1}2
       \Big(\sin\frac{\thet_2}2 \cos\frac{\thet_3}2
           + \cos\frac{\thet_2}2 \sin\frac{\thet_3}2 \Big), \\
     \cos^2\frac{\thet_1}2
 &=& -\cos\frac{\thet_1}2 \cos\frac{\thet_2+\thet_3}2
  =  \cos\frac{\thet_1}2
       \Big(\sin\frac{\thet_2}2\sin\frac{\thet_3}2
           - \cos\frac{\thet_2}2 \cos\frac{\thet_3}2 \Big).
\end{eqnarray*}
dividing by $\prod\sin(\thet_i/2)$, we obtain the third identity
(\ref{iden}). \qed

\section{Discussion}

Our computation of probability distributions in Theorems \ref{theo:int} and \ref{theo:circ} are based on the assumptions that length of a tooth is uniformly distributed on $[0,1]$ and the teeth are situated on the base of the comb with a constant step. However we believe that the answers (i.e., formulas for densities of these distributions) are not sensitive to these details. For example, the answers should be the same if the teeth are randomly uniformly distributed on the base and/or lengths of the teeth are identical independently distributed random values with a continuous distribution density on $[0,b]$, $0<b<\infty$. It would be interesting to obtain a proof of this conjecture.


We have already mentioned that it is interesting to consider analogous problems where base of a random comb is two-dimensional, for example, a disk. Also we would be happy to see dynamical applications of these problems.

\end{document}